\documentclass[reqno]{amsart}
\usepackage{amsmath,amsthm,amssymb,amscd}
\usepackage[pagebackref]{hyperref}

\newcommand{\N}{\mathbb N}
\newcommand{\Z}{\mathbb Z}
\newcommand{\Q}{\mathbb Q}

\newfont{\cyr}{wncyb10}

\newcommand{\rsp}{\raisebox{0em}[2.3ex][1.9ex]{\rule{0em}{2ex} }}

\newcounter{lemmacount}[section]

\newtheorem{thm}[lemmacount]{Theorem}
\newtheorem{prop}[lemmacount]{Proposition}
\newtheorem{lem}[lemmacount]{Lemma}

\title[Higher Descent on Pell Conics]{Higher Descent on Pell Conics. \\
      II. Two Centuries of Missed Opportunities}
\author{Franz Lemmermeyer}
\address{Department of Mathematics,
  Bilkent University,
  06800 Bilkent, Ankara, Turkey}
\email{franz@fen.bilkent.edu.tr}

\begin{document}
\maketitle

\section*{Introduction}

It was already observed by Euler \cite{Euler} that the method
of continued fractions occasionally requires a lot of tedious 
calculations, and even Fermat knew -- as can be seen from the 
examples he chose to challenge the English mathematicians -- a 
few examples with large solutions. To save work, Euler suggested 
a completely different method, which allows to compute even very 
large solutions of certain Pell equations rather easily; its 
drawback was that the method worked only for a specific 
class of equations.  Although Euler's tricks were rediscovered 
on an almost regular basis, nobody really took this approach 
seriously or generalized it to arbitrary Pell equations.

The main goal of \cite{L1} and this article is to discuss certain 
results that have been obtained over the last few centuries and
which will be put into a bigger perspective in \cite{L2}. This is 
opposite to what Dickson aimed at when he wrote his history; in 
\cite[vol II, preface]{Dick} he says
\begin{quote}
What is generally wanted is a full and correct statement of
the facts, not an historians personal explanation of those facts.
\end{quote}
Dickson's books have occasionally been criticised for putting
trivial results next to important ideas, that is, for not 
separating the wheat from the chaff; observe, however,
that a history concentrating only on important ideas would 
hardly have mentioned most of the references on Pell's 
equation that are important for us. This is not because
Dickson failed to see their importance, but because without
the framework of a general theory they were hardly more than 
``mildly amusing'', as van der Poorten \cite{vdP} puts it in 
his review of \cite{Art}, the most recent paper containing 
examples of what we will call second $2$-descents on certain 
Pell conics.

In the follow-up \cite{L2} to this article, I will explain the
theory of the first $2$-descent, study parts of the Selmer
and Tate-Shafarevich groups attached to Pell conics, and
interpret the results discussed here from this modern
point of view. 

\section{Euler}\label{SE}
\subsection{The Content of Euler's Article}
In \cite{Euler}, Euler remarked that the method for solving 
Pell's equation $x^2 - dy^2 = 1$ he has given in \cite{Eul4}
(essentially equivalent to the CFM) is very powerful, but 
that there are certain values of $d$ for which even the CFM 
produces the solution only after tedious calculations; he 
mentions the example $d = 61$, where the smallest positive 
solution is given by $x = 1766319049$ and $y = 226153980$.

Euler then goes on to describe how solutions of certain
auxiliary equations lead to solutions of the Pell equation.
His first equation is
\begin{equation}\label{EE1}
q^2 - ap^2 = -1.
\end{equation}

He writes
\begin{quote}
Problema 1. Si fuerit $app-1 = qq$, invenire numeros $x$ et $y$,
ut fiat $axx+1 = yy$.\footnote{Problem 1. Given $ap^2-1 = q^2$, 
to find numbers $x$ and $y$ such that $axx+1 = yy$.}
\end{quote}
Euler multiplies $ap^2-1 = q^2$ through by $4q^2$ and adds $1$ to get
$4ap^2q^2 + 1 = 4q^4 + 4q^2 + 1 = (2q^2+1)^2$. He has proved:

\begin{lem}\label{EL1}
If $-1 = q^2 - ap^2$, then $y^2 - ax^2 = 1$ for $x = 2pq$ and $y = 2q^2+1$.
\end{lem}

Now Euler investigates the new equation $ap^2-1 = q^2$ more closely:
\begin{quote}
Problema 2. Investigare numeros $a$, pro quibues fieri potest 
$app - 1 = qq$, hincque ipsos numeros $x$ et $y$ assignare, 
ut fiat $axx+1 = yy$.\footnote{Problem 2. To investigate the 
numbers $a$ for which we can solve $ap^2-1 = q^2$, and then to
assign the numbers $x$ and $y$ satisfying $ax^2 + 1 = y^2$.}
\end{quote}

Euler writes $ap^2 = q^2+1$ in the form $a = \frac{q^2+1}{p^2}$
and observes that he has to find $p$ and $q$ in such a way that
this fraction becomes an integer. Since $a$ and $p$ (and therefore
also $p^2$) are sums of two squares, there exist integers 
$b, c, f, g$ such that $p^2 = b^2 + c^2$ and 
$q^2 + 1 = (b^2+c^2)(f^2+g^2)$, so in particular $a = f^2 + g^2$.
Comparing both sides he deduces that, for an appropriate choice
of these numbers, we must have $q = bf+cg$ and $\pm 1 = bg - cf$.

Now Euler assigns values to $b, c, f, g$ in such a way that 
$\pm 1 = bg - cf$ (which can be done in infinitely many ways),
and then computes $a$, $p$, $q = bf+cg$ and finally $x = 2pq$ 
and $y = 2q^2+1$. Actually, he starts by fixing $p = 5$,
which leads to $b = 3$, $c = 4$ in view of $p^2 = 25 = 3^2 + 4^2$.
Then he computes the following table:

\begin{table}[ht!]
$$\begin{array}{r|r|r|r|r|r|r}
  f &  1 &   2 &    4 &    5 &    7 &    8 \\
  g &  1 &   3 &    5 &    7 &    9 &   11 \\ \hline 
  a &  2 &  13 &   41 &   74 &  130 &  185 \\ 
  q &  7 &  18 &   32 &   43 &   57 &   68 \\
  x & 70 & 180 &  320 &  430 &  570 &  680 \\
  y & 99 & 649 & 2049 & 3699 & 6499 & 9249
  \end{array} $$
\caption{}\label{T1}
\end{table}
He also remarks that the cases in this table are not very difficult,
and gives the additional examples $p = 13$, $p=17$, and $p = 25$.

For ease of reference, let us collect Euler's result in the
following

\begin{prop}\label{EP1}
The equation $-1 = q^2 - ap^2$ is solvable if and only if there 
exist $f, g \in \N$ with $a = f^2 + g^2$ such that there are 
$b, c \in \N$ with $p^2 = b^2 + c^2$, $bg - cf = \pm 1$, and 
$q = bf + cg$.
\end{prop}

Now Euler turns to the next equation, namely 
\begin{equation}\label{EE2}
q^2 - ap^2 = -2.
\end{equation}
He first observes (\cite[Problema 4]{Euler})

\begin{lem}\label{EL2}
If $-2 = q^2 - ap^2$, then $y^2 - ax^2 = 1$ for $x = pq$ and $y = q^2+1$.
\end{lem}

As above he then deduces from $ap^2 = q^2 + 2$ that 
$a = f^2 + 2g^2$ for some $f, g \in \N$, and that 
$p^2 = b^2 + 2c^2$, and then concludes that $cf - bg = \pm 1$. 

\begin{prop}\label{EP2}
The equation $-2 = q^2 - ap^2$ is solvable if and only if there 
exist $f, g \in \N$ with $a = f^2 + 2g^2$ such that there are 
$b, c \in \N$ with $p^2 = b^2 + 2c^2$, $bg - cf = \pm 1$, and 
$q = bf + 2cg$.
\end{prop}

Euler's first example is $p = 3$, which in view of 
$3^2 = 1^2 + 2 \cdot 2^2$ implies $b=1$ and $c = 2$. 
Euler solves $2f-g = \pm 1$, giving $a = f^2 + 2g^2$,
$q = f+4g$, hence the solutions $x = 3q$ and $y = q^2+1$
of the Pell equation $1 = y^2 - ax^2$:

\begin{table}[ht!]
$$\begin{array}{r|r|r|r|r|r|r|r}
  f &  1 &   1 &    2 &    2 &    3 &    3 &    4 \\
  g &  1 &   3 &    3 &    5 &    5 &    7 &    7 \\ \hline 
  a &  3 &  19 &   22 &   54 &   59 &  107 &  114 \\ 
  q &  5 &  13 &   14 &   22 &   23 &   31 &   32 \\
  x & 15 &  39 &   42 &   66 &   69 &   93 &   96 \\
  y & 26 & 170 &  197 &  485 &  530 &  962 & 1025
  \end{array} $$
\caption{}\label{T2}
\end{table}

In addition, Euler discusses the examples $p = 9, 11, 17, 19$
and then (Problema 5) goes on to investigate the equation 
\begin{equation}\label{EE3}
q^2 - ap^2 = +2.
\end{equation}

\begin{lem}\label{EL3}
If $2 = q^2 - ap^2$, then $y^2 - ax^2 = 1$ for $x = pq$ and $y = q^2-1$.
\end{lem}

His main result in this case is
 
\begin{prop}\label{EP3}
The equation $2 = q^2 - ap^2$ is solvable if and only if there 
exist $f, g \in \N$ with $a = f^2 - 2g^2$ such that there are 
$b, c \in \N$ with $p^2 = b^2 - 2c^2$, $bg - cf = \pm 1$, and 
$q = bf - 2cg$.
\end{prop}

As examples, Euler treats the cases $p = 7, 17, 23$.

In Problema 7, Euler studies the equation $app + 4 = qq$:

\begin{lem}\label{EL4}
If $4 = q^2 - ap^2$, and if $p$ and $q$ are odd, then $y^2 - ax^2 = 1$ 
for the integers $x = p\,\frac{q^2-1}2$ and $y = q\,\frac{q^2-3}2$.
\end{lem}

This case will not be of interest to us, so let us go right to
Euler's final case, the equation 
\begin{equation}\label{EE4}
s^2 - ar^2 = -4.
\end{equation}

\begin{lem}\label{EL5}
If $-4 = s^2 - ar^2$, then $y^2 - ax^2 = 1$ for 
$x = p\,\frac{q^2-1}2$ and $y = q\,\frac{q^2-3}2$,
where $p = rs$ and $q = s^2$.
\end{lem}

The results of Euler's problema 10 are collected in the following

\begin{prop}\label{EP5}
The equation $-4 = q^2 - ap^2$ is solvable if and only if there 
exist $f, g \in \N$ with $a = f^2 - 2g^2$ such that there are 
$b, c \in \N$ with $p^2 = b^2 - 2c^2$, $bg - cf = \pm 1$, and 
$q = bf - 2cg$.
\end{prop}

Euler's fourth example concerns $a = 109 = 10^2 + 3^2$, 
\begin{quote}
qui methodo vulgari molestissimos calculos requirit,\footnote{for which
the usual method requires the most tedious calculations.}
\end{quote}
Here he takes $25^2 = 625 = 24^2 + 7^2$, and finds $f = 10$, $g = 3$,
$b = 24$ and $c = 7$, hence $s = 261$, giving $p = 6525$ and
$q = 261^2 + 2 = 68123$ and finally
$$\begin{array}{rcrcr} \medskip
   x & = &  6525 \big(\frac{68123^2-1}{2}\big) & = &  15140424455100, \\
   y & = & 68123 \big(\frac{68123^2-3}{2}\big) & = & 158070671986249.
\end{array} $$

\subsection{Interpretation with Continued Fractions}
In the article discussed above, Euler presents solutions of 
the Pell equation for a variety of discriminants. 
The solutions given in Table \ref{T1} correspond to two 
families, namely $d = (3k-1)^2 + (4k-1)^2$ and 
$d = (3k+1)^2 + (4k+1)^2$. Developing $\sqrt{d}$ into
continued fractions (which Euler did not do) we find:

$$ \begin{array}{c|c}
    d & \sqrt{d} \\ \hline
\rsp   (3k-1)^2 + (4k-1)^2 = 5^2k^2 - 14k + 2 &
         [5k-2, \overline{1, 1, 1, 1, 10k-4}] \\ 
\rsp   (3k+1)^2 + (4k+1)^2 = 25k^2 + 14k + 2 & 
         [5k+1, \overline{2, 2, 10k+2}] \\
\rsp   (5k-2)^2 + (12k-5)^2 = 13^2k^2 - 140k + 29 & 
         [13k-6, \overline{1,1,1,1,1,1,26k-12}]
\end{array} $$
The last line comes from Euler's second table, which we did
not reproduce here.

In general, assume that $(r,s,t)$ is a primitive Pythagorean 
triple with $s$ even. Then there exist $m, n$ such that 
$r = m^2 - n^2$, $s = 2mn$ and $t = m^2 + n^2$. Euler has to 
solve the linear equation $(m^2 - n^2)f - 2mng = \pm 1$.

Let us look at families that admit $n = 1$ as a solution.
Then we have to consider $(m^2-1)f - 2mg = \pm 1$; if we put 
$m = 2k$, we find $(4k^2-1)f - 4kg = \pm 1$. The solutions 
of this equation are 
$f = 1 + 4ku$, $g = k+(4k^2-1)u$ and
$f = -1 + 4ku$, $g = -k + (4k^2-1)u$ for $u = 0, 1, 2, \ldots$.
As above we find 
$$ \begin{array}{c|cl}
   d & \sqrt{d} & \\ \hline
\rsp  {}[(4k^2+1)u+k]^2 + 4ku + 1 
           & [e,\overline{2k,2k,2e}] & e = (4k^2+1)u+k \\
\rsp  {}[(4k^2+1)u-k]^2 - 4ku + 1 
           & [e,\overline{1,2k+1,2k+1,1,2e}\,] & e = (4k^2+1)u - k - 1 
   \end{array} $$

\section{Hart}

\subsection{Hart's Article}
In a short note \cite{Hart}, Hart studied the solvability of the
diophantine equation $X^2 - dY^2 = -1$; let us quote the first 
half of \cite{Hart}:
\begin{quote}
To find general values of $x$ and $y$ to solve the problem
$$ x^2 - Ay^2 = -1. $$
Let $A$ be a non-quadrate number $=$ the sum of two squares
$= r^2 + s^2$, then we shall have
$$ x^2 - (r^2 + s^2)y^2 = -1, $$
and by transposition
$$ x^2 - r^2y^2 = s^2y^2 - 1 = (sy-1)(sy+1); $$
whence
$$ x^2 = r^2y^2 + (sy-1)(sy+1) = \square; \quad . \dot{\phantom{p}} . \
   \ x = [r^2y^2 + (sy-1)(sy+1)]^{1/2}.$$
Let $r^2y^2 + (sy-1)(sy+1) = [ry - (m \div n)(sy-1)]^2$. Reducing this
we get
$$ y = \frac{m^2 + n^2}{sm^2 - 2nrm - sn^2}. $$
Here, in order to have $y$ integral, put $sm^2 - 2nrm - sn^2 = \pm 1$;
whence, transposing and dividing by $s$, we have
$$ m^2 - \frac{2nr}s m = \frac{sn^2 \pm 1}{s}, $$
and by quadratics
$$ m = \frac{rn \pm \sqrt{(r^2+s^2)n^2 \pm s}}{s},
    \quad . \dot{\phantom{p}} . \ \ y = \pm (m^2 + n^2). $$
In the general value of $m$, $r$ and $s$ may represent any numbers
one of which is even, and the other any odd number except $1$, and
$n$ can be found by trial, or, if large, by the solution of the
formula $P^2 - (r^2+s^2)n^2 = \pm s$.
\end{quote}

Let us now give a modern interpretation of Hart's idea. 
Assume we want to solve 
\begin{equation}\label{EH1}
x^2 - Ay^2 = -1
\end{equation}
for some squarefree integer $A$. It was already known to
Brahmagupta (see Whitford \cite{Whit}) that this implies that 
$A$ is the sum of two squares. Thus $A = r^2 + s^2$ for integers 
$r, s$, and we have to solve $x^2 = r^2y^2 + (sy-1)(sy+1)$. 
Now parametrize this conic using the rational point 
$P = (x,y) = (\frac{r}s,\frac1s)$: the equation $x - ry = t(sy-1)$ 
describes lines through $P$ with slope $\frac1{r+st}$; if we
pick $t = \frac{m}{n}$ rational, such a line will intersect
the conic in another rational point, and a simple calculation
shows that its $y$-coordinate is 
$y = \frac{m^2 + n^2}{sm^2 - 2nrm - sn^2}$. We want values 
$m$ and $n$ for which $y$ is integral; thus we are led to 
consider 
\begin{equation}\label{EH2}
    sm^2 - 2nrm - sn^2 = \pm 1.
\end{equation}
Any integral solution of (\ref{EH2}) will give an integral 
solution of (\ref{EH1}).
Euler and Lagrange have shown how to reduce (\ref{EH2}) to a
Pell equation: interpreting (\ref{EH2}) as a quadratic polynomial 
in $m$, a necessary condition for the existence of an integral 
solution is that the discriminant $4r^2n^2 + 4s(sn^2 \mp 1)$ 
be a square, i.e., that one of the equations 
\begin{equation}\label{EH3}
        P^2 - An^2 = \pm s
\end{equation}
have an integral solution. This shows

\begin{prop}
Let $d = r^2 + s^2$ be a sum of two squares; if the 
equation $sm^2 - 2rmn - sn^2 = \pm 1$ has an integral
solution $(m,n)$, then so does $x^2 - dy^2 = -1$; in fact,
we can put $y = m^2 + n^2$.
\end{prop}

Thus Hart has proved that the solvability of (\ref{EH2}) for
some choice of $r, s$ implies the solvability of (\ref{EH1}); 
it does not follow (at least not directly) from Hart's proof 
that the condition is also necessary. The necessity actually 
was proved by Euler: Hart's result is nothing but a special 
case of Proposition\ref{EP1}. In fact, if (\ref{EH1}) is
solvable, then by Euler there exist $r,s$ with $A = r^2 + s^2$ 
such that $y^2 = b^2 + c^2$ and $bs - cr = \pm 1$. Since 
$(b,c,y)$ is a Pythagorean triple, we can write $b = m^2 - n^2$, 
$c = 2mn$ and $y = m^2 + n^2$; plugging the values of $b$ and 
$c$ into $bs - cr = \pm 1$ then yields (\ref{EH2}). Thus we 
have shown that combining the results by Euler and Hart gives

\begin{prop}
The equation $x^2 - Ay^2 = -1$ is solvable if and only if
there exist $r, s \in \N$ with $A = r^2 + s^2$ such that 
the diophantine equation (\ref{EH2}) has an integral solution.
\end{prop}

\subsection{The Negative Pell Equation}
Hart's equations were mentioned by Dickson \cite{Dick} and 
rediscovered by Sansone \cite{San,San2} and Epstein \cite{Eps};
these two authors used the fact that solvability of the negative
Pell equation implies the solvability of $x^2 - dy^2 = a$
(Hart's equation (\ref{EH3})),
where $d = a^2 + b^2$ and $a$ is odd. Reducing this equation 
modulo primes $p$ dividing $d$ shows that $(a/p) = +1$:

\begin{prop}\label{PSE}
Let $d$ be a squarefree integer. If $x^2 - dy^2 = -1$ is solvable, 
then $d = a^2 + b^2$ with $a$ odd, and $(a/p) = +1$ for all
primes $p \mid d$.
\end{prop}   

Escott \cite{Esc} had already proved that if $x^2 - Dy^2 = -1$ is 
solvable in integers and $D = a^2 + b^2$, then $(a/D) = 1$ or 
$(b/D) = 1$. He also showed that the condition is not sufficient 
by using the example $2306 = 41^2 + 25^2$. Recently, Proposition
\ref{PSE} was rediscovered by Khessami Pilerud \cite{KP}.

Grytczuk, Luca, \& Wojtowicz \cite{GLW} proved the 
following result, which is just a slightly reformulated
version of Euler's Proposition \ref{EP1}:

\begin{prop}
The equation $s^2 - dr^2 = -1$ has an integral 
solution if and only if there is some odd integer $B \in \Z$ 
with $d = A^2 + B^2$ and a Pythagorean triple $(a,b,c)$ such 
that $aA - bB = \pm 1$. 
\end{prop}

\section{Sylvester}

Sylvester \cite{Syl} proved the following

\begin{prop}\label{PSy}
Let $A = 2f^2+g^2$ be a prime with $f$ odd; then the 
equation $fy^2 + 2gxy - 2fx^2 = \pm 1$ is solvable in
integers. 
\end{prop}

For the proof, let $(u,v)$ be the minimal solution of the
Pell equation $u^2 - Av^2 = 1$. Playing the usual game,
Sylvester arrives at $p^2 - Aq^2 = 1$ or  $p^2 - Aq^2 = -2$,
the other signs being excluded because of $A \equiv 3 \bmod 8$.
The minimality of the solution implies that we have $p^2 - Aq^2 = -2$,
and using unique factorization in $\Z[\sqrt{-2}\,]$ we get
$p + \sqrt{-2} = (g+f\sqrt{-2}\,)(y+x\sqrt{-2}\,)^2$. Comparing
the imaginary parts shows that  $fy^2 + 2gxy - 2fx^2 = \pm 1$.

Sylvester's result can be derived from Euler's Proposition
\ref{EP2} in the same way we deduced Hart's result from 
Proposition \ref{EP1}; all we have to do is replace Pythagorean 
triples by solutions of the equation $x^2 + 2y^2 = z^2$.

\section{G\"unther}

In \cite{Gu82}, S. G\"unther discussed the comments of 
Theon Smyrnaeus on Plato's work and concluded that he must
have been familiar with the Pell equation 
\begin{equation}\label{EGu1}
  2x^2-1\ =\ y^2.
\end{equation}
He then tries to reconstruct a possible approach to the
solution of this equation. 

G\"unther suggests writing (\ref{EGu1}) as $x^2 - 1 = y^2 - x^2$,
substitutes $x+1 = \frac{p}q (y+x)$ and $x-1 = \frac{q}p(y-x)$, 
then sets 
\begin{equation}\label{EGu2}
   p^2+2pq-q^2 = z,
\end{equation}
solves this equation for $p$, and then deduces 
$$  x = \frac{ 4q^2+z\mp 2q\sqrt{2q^2+z}}{z}, \quad
    y = \frac{-4q^2-z\pm 4q\sqrt{2q^2+z}}{z}. $$
Thus any solution of (\ref{EGu2}) with $z = \pm 1$ will lead to
an integral solution of the negative Pell equation (\ref{EGu1}).

The same method, he remarks, works for the more general equation
$$(a^2+b^2)x^2-1\ =\ y^2,$$
which can be written in the form
$(ax-1)(ax+1) = (y-bx)(y+bx)$.
Substituting $ax+1 = \frac{p}q (y+bx)$ and $ax-1 = \frac{q}p(y-bx)$ 
leads to equations equivalent to Hart's (\ref{EH2}) and (\ref{EH3}),
but the formulas derived by G\"unther are incorrect, and his final
conclusion is unclear.

\section{G\'erardin}

While extending existing tables of solutions of the Pell equation
by Legendre, Bickmore, and Whitford, A.~G\'erardin \cite{Ger}
complained about the tedious work necessary when using the 
theory of continued fractions:
\begin{quote}
La recherche pratique de la solution minima \'etait faite jusqu'\`a 
pr\'esent sur les fractions continues, ce qui demande en g\'en\'eral
beaucoup de soins et de temps.\footnote{The practical search for the
minimal solution was made up until now via continued fractions, which
demands in general a lot of care and time.}
\end{quote}
He then presents, by giving a few examples, a `new method' for
solving Pell equations; he uses the auxiliary equations
\begin{equation}\label{EG0}
    x^2 - Ay^2 = \pm 4, \pm 2, -1 
\end{equation}
and remarks that from their solutions one can pass easily to 
the solution of the corresponding Pell equation.

For $A = 941 = 29^2 + 10^2$ he solves the pair of equations
\begin{align*}
  31^2 - 941 \cdot 1^2 & =  +2 \cdot 10, \\
 184^2 - 941 \cdot 6^2 & =  -2 \cdot 10, 
\end{align*}
and then constructs a solution of $x^2 - 941y^2 = -4$ by
setting $y = 1^2 + 6^2 = 37$, giving $x = 1135$.

The general cases are given by the following formulas:
\begin{enumerate}
\item[a)] $z^2 - At^2 = -1$, $A = m^2 + n^2$, $t = \alpha^2 + \beta^2$
      \begin{align}
      \label{EG1a} (m\alpha-n\beta)^2 - A \beta^2  & = \pm m, \\ 
      \label{EG1b} (m\beta+n\alpha)^2 - A \alpha^2 & = \mp m.
      \end{align}
\item[b)] $z^2 - At^2 = +2$, $A = m^2 - 2n^2$, $t = \alpha^2 - 2\beta^2$
      \begin{align}
      \label{EG2a} (n\alpha-m\beta)^2    - A \beta^2  & = \pm n, \\
      \label{EG2b} (n\beta-m\alpha)^2 - A \alpha^2 & = \pm 2n.
      \end{align}
\item[c)] $z^2 - At^2 = -4$, $A = m^2 + n^2$, $t = \alpha^2 + \beta^2$
      \begin{align}
      \label{EG3a} (n\alpha - m\beta)^2 - A \beta^2 & = \pm 2n, \\
      \label{EG3b} (n\beta + m\alpha)^2 - A \alpha^2 & = \mp 2n.
      \end{align}
\end{enumerate}
G\'erardin also remarks that the case $z^2 - At^2 = -2$ can be 
treated similarly. 

Thus in case a), G\'erardin writes $A = m^2 + n^2$ and then
tries to solve the two equations $r^2 - As^2 = m$ and
$t^2 - Au^2 = -m$; putting $t = \alpha^2 + \beta^2$ then gives a
solution of $z^2 - At^2 = -1$, from which a solution to
the Pell equation $X^2 - AY^2 = 1$ is easily derived.
Observe the similarity with the result of Hart \cite{Hart} 
discussed above.

G\'erardin does not give any proofs, but his claims are
easily verified. Let us first consider equation (\ref{EG1a}). 
Plugging in $A = m^2 + n^2$ and simplifying we get
\begin{equation}\label{EG4a}
 m\alpha^2 - 2n\alpha\beta  - m \beta^2 = \pm 1. 
\end{equation}  
Observe that taking (\ref{EG1b}) would lead to the very same equation;
in particular, (\ref{EG1a}) and (\ref{EG1b}) are equivalent.
Also note that $m$ must be odd for (\ref{EG4a}) to be solvable.
A simple calculation now shows that 
\begin{align}
\notag 
(n\alpha^2 - 2m\alpha\beta - n\beta^2)^2 - A(\alpha^2+\beta^2)^2 & =  \\
\label{EG4b}  - (m\alpha^2 - 2n\alpha\beta  - m \beta^2)^2       & =  -1. \\ 
\intertext{Similarly, in cases b) and c) we get the equations}
\label{EG5a} n\alpha^2 - 2m\alpha\beta + 2n\beta^2 & =  \pm 1 \\
\label{EG6a} n\alpha^2 - 2m\alpha\beta -  n\beta^2 & =  \pm 2 \\
\intertext{(it is easy to see that, in (\ref{EG6a}), 
            $m$ must be odd) as well as}
\notag (m\alpha^2 - 4n\alpha\beta + 2m\beta^2)^2 
                        - A(\alpha^2-2\beta^2)^2 & =  \\
\label{EG5b} 2(n\alpha^2 - 2m\alpha\beta + 2n\beta^2)^2 & =  +2 \\
\notag  (m\alpha^2 - 2n\alpha\beta + m\beta^2)^2 
                              - A(\alpha^2+ \beta^2)^2 & =  \\ 
\label{EG6b} -(n\alpha^2 - 2m\alpha\beta - n\beta^2)^2 & =  -4. 
\end{align}

Comparing (\ref{EG0}) with Euler's equations (\ref{EE1}) -- (\ref{EE4}) 
from Section \ref{SE} it will come as no surprise that G\'erardin's 
formulas are an easy consequence of Euler's results. 

\section{Hardy \& Williams, Bapoungu\'e, Arteha}

\subsection{Hardy \& Williams}
In \cite{HaWi}, K.~Hardy \& K.~Williams investigate the
solvability of the diophantine equation
\begin{equation}\label{EHW}
  fx^2 - 2gxy - fy^2 = 1.
\end{equation}
Their main result is

\begin{prop}\label{PHW}
The equation $-1 = q^2 - ap^2$ is solvable if and only if there 
exist $f, g \in \N$ with $a = f^2 + g^2$ such that there are 
$x, y \in \Z$ with (\ref{EHW}). In this case, the pair $(f,g)$
with $f$ odd is unique. 
\end{prop}

Apart from the uniqueness assertion, this result is an almost 
trivial consequence of Euler's Proposition \ref{EP1}: we start 
with the observation that $(b,c,p)$ is a Pythagorean triple. 
Now $-1 = q^2 - ap^2$ implies that $a$ and $p$ are odd. Assuming 
that $c$ is odd, we find from $bg - cf = \pm 1$ that $f$ is also 
odd. Changing the signs of $f, g$ if necessary we may assume that 
$bg - cf = -1$. The parametrization of Pythagorean triples shows 
that $b = 2xy$ and $c = x^2 - y^2$; the equation $bg - cf = -1$ 
then becomes $f(x^2-y^2) - 2gxy = 1$. This proves Proposition 
\ref{PHW} except for the uniqueness part.

The claim that there is essentially only one such pair $(f,g)$
is an important contribution: as we will see later, it should
be seen as (part of) an analogue of Dirichlet's Theorem 
\cite[Thm. 3.3]{L1}. 

\subsection{Bapoungu\'e}
Inspired by the work of Hardy \& K.~Williams \cite{HaWi}, Bapoungu\'e 
\cite{Bap,Bap1,Bap2,Bap3,Bap4} started investigating the solvability of 
the diophantine equation
\begin{equation}\label{EB1}
ax^2 + 2bxy - kay^2 = \pm 1
\end{equation}
for values of $k$ for which $\Q(\sqrt{-k}\,)$ has class number $1$.
The identity
$$ (-bx^2 + 2akxy + kby^2)^2 - (b^2+ka^2)(x^2 + ky^2)^2
     = -k(ax^2 + 2bxy - kay^2)^2 $$
yields, upon substituting a solution of (\ref{EB1}) for $(x,y)$,
the equation
\begin{equation}\label{EBk}
(-bx^2 + 2akxy + kby^2)^2 - \delta(x^2 + ky^2)^2 = -k.
\end{equation}
This shows (\cite[Thm. 2]{Bap1})

\begin{thm}
If (\ref{EB1}) has an integral solution, then so does (\ref{EBk}).
\end{thm}

Multiplying (\ref{EB1}) through by $a$ and completing the square we get
\begin{equation}\label{EB2}
 (ax+by)^2 - \delta y^2 = \pm a, \quad \text{where}\ \delta = ka^2+b^2.
\end{equation}
Thus the solvability of (\ref{EB1}) implies the solvability of the 
Pell equation $X^2 - \delta Y^2 = a$, where solvability denotes 
solvability in integers. 

Similarly, multiplying (\ref{EB1}) through by $-ka$ and completing 
the square we get
\begin{equation}\label{EB3}
(kay - bx)^2 - \delta x^2 = \mp ka.
\end{equation}
The special cases $k = 1, 2$ of these equations go back to Euler,
Hart, and G\'erardin. 

The main result of Bapoungu\'e's thesis \cite{Bap} is

\begin{thm}
Let $k \in \{2, 3, 7, 11, 19, 43, 67, 163\}$ (this implies that
$\Q(\sqrt{-k}\,)$ has class number $1$). If $a, b$ are positive
integer with $a$ odd such that $p = ka^2 + b^2$ is prime, then
(\ref{EB1}) is solvable if and only if (\ref{EBk}) is solvable.
\end{thm}

A similar result holds for $k = 1$. The case $k = 2$ is 
Sylvester's Proposition \ref{PSy}.
In \cite{Bap3}, the following result is proved:

\begin{thm}
Let $k$ be as above. Among all pairs $(a,b)$ of coprime natural
numbers with $a$ odd and $d = ka^2 + b^2$, there is exactly one
pair for which (\ref{EB1}) is solvable.
\end{thm}

This generalizes Theorem \ref{PHW}.
 
\subsection{Arteha}

The last rediscovery of the method of Euler-Hart-G\'erardin
so far is due to Arteha \cite{Art}; his results are

\begin{prop}
Consider the Pell equation
\begin{equation}\label{EA1} x^2 - dy^2 = 1 \end{equation}
for primes $d$.
\begin{enumerate}
\item If $d \equiv 1 \bmod 4$, and write $d = a^2 + b^2$
      with $a$ odd. Then the minimal positive solution of Pell's
      equation (\ref{EA1}) is given by 
      $$ y = 2|2amn + b(m^2-n^2)|(m^2 + n^2), $$
      where $m$ and $n$ satisfy
      $$ a(m^2 - n^2) - 2bmn = \pm 1.$$
\item If $d \equiv 3 \bmod 8$, write $d = a^2 + 2b^2$. Then the 
      minimal positive solution of Pell's equation (\ref{EA1}) 
      is given by 
      $$ y = \big|4bmn + a|m^2 - 2n^2|\big|(m^2 + 2n^2), $$
      where $m$ and $n$ satisfy
      $$ b|m^2 - 2n^2| - 2amn = \pm 1.$$
\item If $d \equiv 7 \bmod 8$, write $d = a^2 - 2b^2$. Then the 
      minimal positive solution of Pell's equation (\ref{EA1}) 
      is given by 
      $$ y = |(a(m^2 + 2n^2) - 4bmn)(m^2-2n^2)|, $$
      where $m$ and $n$ satisfy
      $$ 2amn - b(m^2 +2n^2) = \pm 1.$$
\end{enumerate}
\end{prop}

\section{Summary}
Starting with Euler in \cite{Euler}, many authors have come
up with essentially the same idea: solving the Pell equation
\begin{equation}\label{ES1} X^2 - dY^2 = 1 \end{equation}
becomes easier by looking at certain auxiliary equations. 

The first step is writing down Legendre's equations
\begin{equation}\label{ES2}
    rx^2 - sy^2 = 1, 2 \quad \text{for} \ d = rs.
\end{equation}
There is one nontrivial equation among these with a solution, 
and the smallest solution will have about half as many digits 
as the smallest solution of (\ref{ES1}). 

The second step is to look at equations whose solutions give
rise to a solution of one of the equations (\ref{ES2}); but this
second step has never been completed in full generality. 
What we have are specific equations applicable only in 
special situations; these were first discovered by Euler,
rediscovered by Hart, Sylvester, G\"unther, G\'erardin, 
Hardy \& Williams, and Arteha, and slightly generalized 
by Bapoungu\'e.


\begin{thebibliography}{99}

\bibitem[Art2002]{Art} S.N. Arteha,
{\em Method of hidden parameters and Pell's equation},
JPJ Algebra Number Theory Appl. {\bf 2} (2002), 21--46; cf. p.
%

\bibitem[Bap1989]{Bap} L. Bapoungu\'e, 
{\em Sur la r\'esolubilit\'e de l'\'equation
     $ax^2 + 2bxy - kay^2 = \pm 1$}, 
Th\`ese Univ. Caen, 1989; see also
C. R. Acad. Sci. Paris {\bf 309} (1989), 235--238; cf. p.  
%

\bibitem[Bap1998]{Bap1} L. Bapoungu\'e, 
{\em Un crit\`ere de r\'esolution pour l'\'equation diophantienne
     $ax^2 + 2bxy - kay^2 = \pm 1$},
Expos. Math. {\bf 16} (1998), 249--262; cf. p.  
%

\bibitem[Bap2000a]{Bap2} L. Bapoungu\'e, 
{\em Sur la r\'esolubilit\'e de l'\'equation
     $ax^2 + 2bxy - 8ay^2 = \pm 1$}, 
IMHOTEP, J. Afr. Math. Pures Appl. {\bf 3} (2000), 97--111; cf. p.  
%

\bibitem[Bap2000b]{Bap3} L. Bapoungu\'e, 
{\em Sur les solutions g\'enerales de l'\'equation
     diophantienne $ax^2 + 2bxy - kay^2 = \pm 1$}, 
Expos. Math. {\bf 18} (2000), 165--175; cf. p.  
%

\bibitem[Bap2002]{Bap4} L. Bapoungu\'e, 
{\em The diophantine equation $ax^2+2bxy - 4ay^ 2 = \pm 1$},
Intern. J. Math. Math. Sci. {\bf 35} (2003), 2241--2253
%

\bibitem[Dic1920]{Dick} L.E. Dickson,
{\em History of the Theory of Numbers},
vol I (1920); vol II (1920); vol III (1923);
Chelsea reprint 1952; cf. p. 
%

\bibitem[Eps1934]{Eps} P. Epstein,
{\em Zur Aufl\"osbarkeit der Gleichung $x^2-Dy^2=-1$},
J. Reine Angew. Math. {\bf 171} (1934), 243--252; cf. p.   
%

\bibitem[Esc1905]{Esc} E.B. Escott,
{\em Solution de l'\'equation $x^2 - Dy^2 = -1$}, 
L'Interm\'ed Math. {\bf 12} (1905), 53; cf. p.  
%

\bibitem[Eul1765]{Eul4} L. Euler,
{\em De usu novi algorithmi in problemate Pelliano solvendo},
Novi Acad. Sci. Petropol. {\bf 11} (1765) 1767, 28--66;
Opera Omnia I-3, 73--111; cf. p.  
%

\bibitem[Eul1773]{Euler} L. Euler,
{\em Nova subsidia pro resolutione formulae $axx + 1 = yy$},
Sept. 23, 1773; 
Opusc. anal. {\bf 1} (1783), 310; 
Comm. Arith. Coll. {\bf II}, 35--43;
Opera Omnia I-4, 91--104; cf. p.  
%

\bibitem[Ger1917]{Ger} A. G\'erardin,
{\em Sur l'\`equation $x^2 - Ay^2 = 1$}, 
L'Ens. math. {\bf 19} (1917), 316--318;
Sphinx-\OE dipe {\bf 12} June 15, 1917, 1--3; 
cf. p.   
%

\bibitem[GLW2000]{GLW} A. Grytczuk, F. Luca, M. Wojtowicz,
{\em The negative Pell equation and Pythagorean triples},
Proc. Japan Acad. {\bf 76} (2000), 91--94; cf. p.  
%

\bibitem[Gue1882]{Gu82} S. G\"unther, 
{\em Ueber einen Specialfall der Pell'schen Gleichung},
Bl\"atter f\"ur das Bayerische Gymnasial- und Realschulwesen
{\bf 17} (1882), 19--24; cf. p.  
%

\bibitem[HW1986]{HaWi} K. Hardy, K. Williams,
{\em On the solvability of the diophantine equation 
     $dV^2 - 2e VW - dW^2 = 1$},
Pac. J. Math. {\bf 124} (1986), 145--158; cf. p.  
%

\bibitem[Har1878a]{Hart} D.S. Hart, 
{\em Solution of an indeterminate problem},
Analyst {\bf 5} (1878), 118--119; cf. p.  
%

\bibitem[KP1999]{KP} Kh. Khessami Pilerud, 
{\em On the Diophantine equation $x^2 - Ny^2 = -1$}, (Russian)
Vestnik Moskov. Univ. Ser. I Mat. Mekh. (1999), no. 2, 65--67; 
Engl. transl. Moscow Univ. Math. Bull. 
    {\bf 54} (1999), no. 2, 48--49; cf. p.
%

\bibitem[Lem2003a]{L1} F. Lemmermeyer,
{\em Higher Descent on Pell Conics I. From Legendre to Selmer},
preprint 2003; cf. p. 
%

\bibitem[Lem2003b]{L2} F. Lemmermeyer,
{\em Higher Descent on Pell Conics III. The First $2$-Descent},
preprint 2003; cf. p. 
%

\bibitem[vdP2003]{vdP} A. van der Poorten,
{\em Review 2003i:11040}, MathSciNet; cf. p. 
%

\bibitem[San1925a]{San} G.~Sansone,
{\em Sulle equazioni indeterminate delle unit\`a di norma negativa
    dei corpi quadratici reali}, Rend. Acad. d. L. Roma 
(6) {\bf 2} (1925), 479--484; cf. p. 
%

\bibitem[San1925b]{San2}  G.~Sansone,
{\em Ancora sulle equazioni indeterminate delle unit\`a di norma 
  negativa dei corpi quadratici reali}, Rend. Acad. d. L. Roma 
(6) {\bf 2} (1925), 548--554; cf. p. 
%

\bibitem[Syl1881]{Syl} J.J. Sylvester,
{\em Mathematical Question 6243},
Educational Times {\bf 34} (1881), 21--22; cf. p. 
%

\bibitem[Whi1912]{Whit} E.E. Whitford,
{\em The Pell equation},
New York 1912, 193 pp; cf. p. 
%

\end{thebibliography}
\end{document}